\documentclass[12pt]{amsart}
\usepackage{amsfonts}
\usepackage{amsfonts,latexsym,rawfonts,amsmath,amssymb,amsthm}
\usepackage[plainpages=false]{hyperref}
\usepackage{graphicx}

\numberwithin{equation}{section}

\RequirePackage{color}
\theoremstyle{plain}

\newtheorem{theorem}{Theorem}[section]

\newtheorem{definition}{Definition}[section]

\newtheorem{lemma}[theorem]{Lemma}
\newtheorem{proposition}[theorem]{Proposition}
\newtheorem{remark}[theorem]{Remark}
\newtheorem{problem}[theorem]{Problem}

\newcommand{\beq}{\begin{equation}}
\newcommand{\eeq}{\end{equation}}
\newcommand{\beqs}{\begin{eqnarray*}}
\newcommand{\eeqs}{\end{eqnarray*}}
\newcommand{\beqn}{\begin{eqnarray}}
\newcommand{\eeqn}{\end{eqnarray}}
\newcommand{\beqa}{\begin{array}}
\newcommand{\eeqa}{\end{array}}

\def\phi{\varphi}

\begin{document}

\title[A class of nonlinear equations ]
{A class of fully nonlinear equations arising in conformal geometry
}

\author{Li Chen}
\address{Faculty of Mathematics and Statistics, Hubei Key Laboratory of Applied Mathematics, Hubei University,  Wuhan 430062, P.R. China}
\email{chernli@163.com}

\author{Xi Guo}
\address{Faculty of Mathematics and Statistics, Hubei Key Laboratory of Applied Mathematics, Hubei University,  Wuhan 430062, P.R. China}
\email{guoxi@hubu.edu.cn}

\author{Yan He}
\address{Faculty of Mathematics and Statistics, Hubei Key Laboratory of Applied Mathematics, Hubei University,  Wuhan 430062, P.R. China}
\email{helenaig@hotmail.com}

\keywords{Conformal geometry; Fully nonlinear equation; Compact manifolds}

\subjclass[2010]{Primary 35J96, 52A39; Secondary 53A05.}

\thanks{This research was supported by funds from Hubei Provincial Department of Education Key Projects D20171004.}

c

\begin{abstract}
In this paper, we consider a class of fully nonlinear equations
on closed smooth Riemannian manifolds,
which can be viewed as an extension of $\sigma_k$
Yamabe equation. Moreover, we prove local gradient and second
derivative estimates for solutions to these equations and establish
an existence result associated to them.
\end{abstract}

\maketitle

\baselineskip18pt

\parskip3pt

\section{Introduction}

Let $(M, g_0)$  be  a smooth, compact Riemannian manifold of
dimension $n\ge3$.
The $k$-th elementary symmetric polynomial
is denoted by $\sigma_k$:
\[\sigma_k(\lambda)=\sum_{1\le i_1<\cdots<i_k\le n}\lambda_{i_1}\cdots\lambda_{i_k}.\]
Let $V$ be a $(0,2)$ tensor on $(M, g)$.
We define $\sigma_k$- curvature of $V$ by
\begin{equation*}
\sigma_k(V)
\end{equation*}
where   $\sigma_k(V)$ means $\sigma_k$ is applied to the eigenvalues
of $g^{-1} V$.
\begin{definition}
We say $V$ is \textbf{$k$-admissible}
if $\sigma_1(V)>0,\cdots,\sigma_k(V)>0$.
\end{definition}

\begin{remark}
$V$ is \textbf{$k$-admissible} if and only if the vector of
eigenvalues of $g^{-1}V$,  $\lambda=(\lambda_1,\cdots,\lambda_n)$,
lies in $\Gamma_k$, i.e. $V\in\Gamma_k$, where $\Gamma_l$ is an open
convex symmetric cone with vertex at the origin
\begin{equation*}
   \Gamma _{l}=\{\lambda =(\lambda _{1},\lambda _{2},\dots ,\lambda _{n})\in
   \mathbb{R}^{n}|\quad \sigma _{j}\left( \lambda \right) >0,1\leq j\leq l\},
\end{equation*}
$1\leq l\leq n$.
\end{remark}

In our  paper,
 we study the problem of prescribing linear combination of $\sigma_k$-curvature of $V$ on $(M, g)$:
\begin{eqnarray}\label{pro}
\sigma_k(V) +\alpha(x) \sigma_{k-1}(V)=f(x) \quad 3\leq k\leq n.
\end{eqnarray}
where $\alpha(x)\in C^\infty(M)$ and $0<f(x)\in C^\infty(M)$, $V$ is
a $(0,2)$ symmetric tensor on $(M, g)$. It  was introduced in Krylov
\cite{Kry95}. The author considered the case with $\alpha(x) \le 0$
and studied Dirichlet problem of the following degenerate equation
in a $(k-1)$-convex domain $D$ in $\mathbb{R}^n$,
\[\sigma_k(D^2u) =
\sum_{l=0}^{k-1} \alpha_l(x)\sigma_l(D^2u),\] based on his
observation that the above equation is elliptic in admissible
$\Gamma_k$ cone if all the coefficient $\alpha_l(x)\geq 0$.
Recently, in
 \cite{GZ19} the authors studied
 \begin{eqnarray*}\label{p}
\sigma_k(D^2u+uI) +\alpha(x)
\sigma_{k-1}(D^2u+uI)=\sum_{l=0}^{k-2}\alpha_l\sigma_l(D^2u+uI)\quad
\mbox{on} \quad \mathbb{S}^n,
\end{eqnarray*}
which arises in the problem of prescribing convex combination of
area measures \cite{Sch13}.  And they also studied Krylov' equation.
A key new feature they observed is that there is no sign requirement
for the coefficient function $\sigma_{k-1}$. Thus, the proper
admissible set of solutions of the equation is $\Gamma_{k-1}$, not
$\Gamma_k$.

Denote by ${Ric}_{g}$ and ${R}_{g}$ the Ricci and scalar curvatures
of $g$ respectively. If we replace $V$ with the Schouten tensor
$${A}_{g}
   ={\frac{1}{{n-2}}}
      \left(
         {Ric}_{g} - {\frac{{{R}_{g}}}{{2(n-1)}}}g
      \right)$$ and take $\alpha(x)=0$ in (\ref{pro}), it turns out to be
a prescribing $\sigma_k$-scalar curvature equation which
was first introduce by \cite{V00} and \cite{V02}:
\begin{equation}\sigma_k(A_g)=f.\label{eq000}\end{equation}
When $f=const.$, (\ref{eq000}) is just the
$\sigma_k$-Yamabe equation.
In the past two decades, a great deal of mathematical effort has been
devoted to the study of $\sigma_k$-Yamabe problem.
For $k=1$ and $f=const.$, (\ref{eq000}) reduce to classical Yamabe problem.
It has been solved by  \cite{Y60,T68,Au76,S84,E92}.
In 2009, a remarkable result was proved by\cite{BM09,KMS09}
 that the solution set of (\ref{eq000}) for $k=1$ is compact if and only if $n\le24$.
 For $k\ge2,$ Viaclovsky  \cite{V00} established an important result, that the metric $g_0$
 is a critical point of the functional
 $g\mapsto \int_M\sigma_k(A_g)dv_g$
 restricted to the set of unit volume metrics
 if and only if $\sigma_k(A_{g_0})=const.$,
 provided either (i) $k\neq n/2$ and $(M,g_0)$ is locally conformally flat
 or (ii) $k=2$.
 In dimension 4, \cite{CGY02a} prove that
  if the Yamabe constant $\inf_{\tilde{g}\in[g],\atop vol(\tilde{g})=1}\int_M \sigma_1(A_{\tilde{g}})$ and
  $\int_M \sigma_2(A_{g_0})$ are both positive,
  then we can find a conformal metric $\tilde{g}$
  such that $\sigma_2(A_{\tilde{g}})$ is a positive constant.
  Later on,
  Guan-Wang and Li-Li proved
 that on  locally conformally flat manifolds,
 the $\sigma_k$- Yamabe equation
 is solvable for $k<n/2$ in \cite{GW03b} and \cite{LL03}
 respectively
 if $\sigma_i(A_{g_0})>0$ for $1\le i\le k$;
 see also \cite{STW07}.
 In \cite{BV04}, the authors present a conformal
 variational characterization in dimension $n =2k$ of the equation
 $\sigma_k(A_g)=constant$.
  In \cite{HL04}, the authors  showed that a metric $g$  on an $n$-dimensional ($n\ge5$)
 locally conformally flat manifold with volume one has constant sectional curvature
 if and only if $g$  is a critical point of functional
\[\mathcal{F}_2(g)=\int\sigma_2(A_g)dv_g\]
 restricted to the metrics of volume one and  $\mathcal{F}_2(g)>0$.
  In \cite{GV07} Gursky and Vioclovsky solved the prescribing
  curvature problem  for $k>n/2$  in
  affirmation,  provided
  $\sigma_i(A_{g_0})>0$ for $1\le i\le k$ and $(M, g)$ is not conformally
  equivalent to a unit sphere.
  In \cite{HLS08}, the authors proved some rigidity properties in terms of the
  curvature $\sigma_k(A_g)$ on closed locally conformally flat manifolds under
  the assumption that $A_g$ is semi-positive definite for
  $k \ge 3$ and $Ric_g$ is semi-positive definite.
 A remarkable result proved by Trudinger and Wang in \cite{TW10} asserts that
 the existence
 theorem  and compactness of solutions to the
 $\sigma_k-$ Yamabe problem for $k=n/2$
 hold true.
  For more general class of functions, including quotient curvature,
  were first considered by Lichnerowiczin\cite{L58}.
 Then in \cite{GLW10}, the authors consider the problem
 \[\frac{\sigma_2}{\sigma_1}(A_g)=f,\] where $f$ is a given smooth function.
 They proved that,
  let $g_0$ be a positive scalar curvature metric,
  then in its conformal class
  there is a conformal metric with $\sigma_2(A_g) = \kappa\sigma_1(A_g)$,
  for some constant $\kappa$.
For the manifolds with boundary,
Chen showed the variational property of the k-curvature and solved
the k-Yamabe problem on locally conformally manifolds with
umbilic boundary in \cite{Chen09}. Related works include
\cite{Chen07, HSh11}, etc..

For the
manifolds with negative curvature,
there are also many interesting results.
In \cite{GV03a}, the authors proved that every compact
manifold with negative Ricci curvature is conformal to
a metric $g$ with $\det{Ric_g}=const.$
In \cite{LiSh05}, the authors using a parabolic approach
 proved that
there exists a conformal metric $g$
such that $\sigma_k(-Ric_g)=const.$
if $-Ric_{g_0}\in \Gamma_k$.
For the manifolds with totally geodesic boundary,
\cite{ShY14} showed that
for some positive function $f(x)$ and $h<k$
there exists a conformal metric $g$
such that
\begin{equation}
\left\{
\begin{array}
[c]{ll}%
\frac{\sigma_k}{\sigma_h}(-Ric_g)=f(x)\label{qeq}\\
h_g=0
\end{array}
\right.
\end{equation}
where $h_g$ is the mean curvature.
In \cite{Sui17}, it was proved that on $\mathbb{R}^n$,
there exists a conformal metric $g$
such that
\begin{equation}\frac{\sigma_k}{\sigma_h}(-Ric_g)=\psi(x),\label{qeq1}\end{equation}
where $k>h$, $\psi(x)$ is a positive function
and $C_1|x|^{-kl}\le \psi(x)\le C_2|x|^{-kl}$ for $x$
large and $l>2$.

A key issue for the study of $\sigma_k$-Yamabe problem is the a priori estimates.
$C^{1}$ and $C^{2}$ estimates for  the
 equation (\ref{eq000}) and (\ref{qeq})
have been studied extensively.
See
\cite{Chen05, GW03a, GW03b, LL03, STW07, W06} for local interior estimates and \cite{V02}
for global estimates on closed manifolds.
Another interesting problem is to study the fully nonlinear equation
(\ref{eq000}) on a compact Riemannian manifold $(M^{n},g)$ with boundary
$\partial M $. In \cite{G07}, Bo Guan studied the existence problem under the
Dirichlet boundary condition. The pioneering works on the Dirichlet
problems for fully nonlinear elliptic equations are \cite{CNS85, Tr90} etc.. The
Neumann problem for (\ref{eq000}) has been studied by S. Chen \cite{Chen07, Chen09},
Jin-Li-Li \cite{JLL07},  Li-Li \cite{LL06},He-Sheng\cite{HSh13} and Sheng-Yuan\cite{ShY13}, etc..
They obtain local estimates for solutions and establish some existence
results under various conditions.

Motivated by Krylov \cite{Kry95} and Guan-Zhang \cite{GZ19}, it is
natural to study the equation of type \eqref{pro} with $V=A_g$ and
$-Ric_g$ which can be viewed as extensions of prescribing $\sigma_k$
curvature problem  in conformal geometry.  Let us  first  consider $V=-\frac{Ric_g}{n-2}$
and study:
\begin{equation}\label{problem2}
\sigma_k(\frac{-Ric_g}{n-2})+\alpha(x)\sigma_{k-1}(\frac{-Ric_g}{n-2})=f(x), \quad 3\leq k\leq n.
\end{equation}

\begin{theorem}\label{main2} Let $(M, g_0)$
be a smooth closed Riemannian manifold of dimension $n\ge3.$ Assume
that $-Ric_{g_0}$ is $k$-admissible for $k\ge 3$. Then there exists
a conformal metric $g$ satisfies equation (\ref{problem2})

if

either

\textbf{Case(A)}
$0\ge\alpha\in C^\infty(M)$,
and  $0<f(x)\in C^\infty(M)$,

or

\textbf{{Case(B)}} $0>\alpha\in C^\infty(M)$, and  $f(x)=0$.
\end{theorem}

To prove Theorem \ref{main2}, we will establish
the following estimates:

\begin{theorem}
Let $(M, g_0)$ be a smooth closed Riemannian manifold of dimension
$n\ge3.$ \label{nest}
 \textbf{Case(A)} Assume $f\ge\theta>0$. Let $u$ be a solution of (\ref{problem2}), $g=\exp(2u)g_0$ and $-Ric_g$  $\in\Gamma_{k-1}$.
Then there exists constant $C$,
 depending on
 $||\alpha||_{C^2(B_r)},
 ||f||_{C^2(B_r)}$, $\sup_{B_r} |u|$, $\theta,$ $g_0$,
such that
\begin{equation}
\sup_{B_{\frac{r}{2}}} |\nabla u|^2+|\nabla^2 u|\le C.\label{2}
\end{equation}

\textbf{Case(B)} Assume $\alpha<0, f\equiv0$. Let $u$ be a solution
of (\ref{problem2}), $g=\exp(2u)g_0$ and $-Ric_g$ $\in\Gamma_{k}$.
Then there exists constant $C$,
 depending on
 $||\alpha||_{C^2(B_r)}
 $, $\sup_{B_r} |u|$, $g_0$,
such that
\begin{equation}
\sup_{B_{\frac{r}{2}}}|\nabla u|^2+|\nabla^2 u|\le C.\label{}
\end{equation}
\end{theorem}

Then let us  replace $V$ with Schouten tensor $A_{{g}}$ in (\ref{pro}) and consider:
\begin{eqnarray}\label{problem}
\sigma_k( A_{{g}}) +\alpha(x) \sigma_{k-1}(A_{{g}})=f(x), \quad 3\leq k\leq n.
\end{eqnarray}

We will establish the following a priori estimates.
\begin{theorem}\label{pes}
Let $(M, g_0)$
be a smooth closed Riemannian manifold of dimension $n\ge3.$
Assume that Schouten tensor $A_{g_0}$ is $k-1$-admissible for $ k\ge 3$,
and $f(x)\ge\theta>0$.
Let $u$ be a solution of (\ref{problem}),
$g=\exp(-2u)g_0$ and $A_g$  $\in\Gamma_{k-1}$.
Then there exists constant $C$,
 depending on
 $||\alpha||_{C^2(B_r)},
 ||f||_{C^2(B_r)}$, $\sup_{B_r} |u|$, $\theta,$ $g_0$,
such that
\begin{equation}
\sup_{B_{\frac{r}{2}}} |\nabla u|^2+|\nabla^2 u|\le C.\label{3}
\end{equation}
\end{theorem}

Then naturally the following interesting problem would be proposed.
\begin{problem}
Since on the locally conformally flat manifolds, $C^0$ estimates
and the existence result for
$\sigma_k$- Yamabe problem hold true,
is the existence result also true for the equation (\ref{problem})
on the locally conformally flat manifolds?

\end{problem}
The present paper is built up as follows.
In Sect. 2
we start with some background.
We will prove $C^0$ estimates in Sect. 3.
The $C^1$ and $C^2$ estimates are treated in Sect. 4.
Theorem \ref{main2} is proved in Sect. 5.

\vskip30pt

\section{Preliminaries}

{We first recall the following Newton-Maclaurin inequality.}
\begin{lemma}\textit{(\cite{Tr90, LT94})} \label{lemma1}
\textit{\ Let} $\lambda\in\mathbb{R}^n$. \textit{ For } $0\leq
l<k\leq n,$ $r>s\ge0, k\ge r, l\ge s$,

(1)
\[
k(n-l+1)\sigma_{l-1}(\lambda)\sigma_{k}(\lambda)\leq l(n-k+1)\sigma_{l}(\lambda)\sigma_{k-1}(\lambda).
\]

(2)
\[\big[\frac{\sigma_k(\lambda)/C^k_n}{\sigma_l(\lambda)/C^l_n}\big]^\frac{1}{k-l}
\le\big[\frac{\sigma_r(\lambda)/C^r_n}{\sigma_s(\lambda)/C^s_n}\big]^{\frac{1}{r-s}},
\quad \textit{for }  \lambda\in \Gamma_k.\]
\end{lemma}
$\Box$

Let  ${g}=\exp{(2u)}g_0$. Then
\begin{eqnarray*}
&&Ric_{{g}}\\
&=&(n-2)\left(  -\nabla^2 u -\frac{1}{n-2} \triangle u g_0 - |\nabla u|^2 g_0+ du\otimes du +\frac{Ric_{g_0}}{n-2} \right).
\end{eqnarray*}
 Denote
 \[U= \nabla^2 u +\frac{1}{n-2} \triangle u g_0 +|\nabla u|^2 g_0- du\otimes du -t\frac{Ric_{g_0}}{n-2}+\frac{1-t}{n}g_0.\]
 We will use the above conformal translation later.

\textbf{Case(A)}We first assume $f>0$. Consider a family of equations
 \begin{eqnarray}
G(U)&=&\frac{\sigma_k\big(tU +(1-t)trU g_0\big)}
    {\sigma_{k-1}\big(tU +(1-t)trU  g_0\big)}\notag\allowdisplaybreaks\\
    &&-\frac{((1-t)\sigma_k(e)+tf(x))\exp({2ku})}
    {\sigma_{k-1}\big(tU +(1-t)trU  g_0\big)}\notag\allowdisplaybreaks\\
 &=&-t\alpha \exp({2u}),\label{eq99}
\end{eqnarray}
where  $ t\in [0, 1], e=(1,\cdots, 1)$.
Then we will prove that the equations we discussed are elliptic and concave.
\begin{proposition}\textit{(\cite{GZ19})} \label{ellipticconcave}
Let  $\eta=(\eta_1,\cdots,\eta_n), e=(1,\cdots,1)$.
Denote by $tr\eta$ the function $\sum_{i=1}^n\eta_i$. Then
the operator
\begin{eqnarray}
&&G(\eta)\allowdisplaybreaks\notag\\
&:=&          \frac{\sigma_k(t\eta+(1-t)tr\eta \cdot e)}{\sigma_{k-1}(t\eta+(1-t)tr\eta \cdot e)}\allowdisplaybreaks\          -h(x)\frac{1}{\sigma_{k-1}(t\eta+(1-t)tr\eta \cdot e)}
\end{eqnarray}
is elliptic and concave
if
$\eta\in \Gamma_{k-1}$,
$0<h(x)\in C^\infty(M)$.
\end{proposition}

\textbf{Proof.}

(1)
Set $\mu=t\eta+(1-t)tr\eta \cdot e$,
$G_k(\eta)= \frac{\sigma_k(t\eta+(1-t)tr\eta \cdot e)}{\sigma_{k-1}(t\eta+(1-t)tr\eta \cdot e)}$.
By direct calculation, we obtain
\begin{eqnarray}
&&   \sigma_{k-1}(\mu|i)\sigma_{k-1}(\mu)-\sigma_k(\mu)\sigma_{k-2}(\mu|i)\notag\allowdisplaybreaks\\
&=&  \sigma_{k-1}(\mu|i)\big(\sigma_{k-1}(\mu|i)+\mu_i\sigma_{k-2}(\mu|i)\big)
     -\big(\sigma_k(\mu|i)+\mu_i\sigma_{k-1}(\mu|i)\big)\sigma_{k-2}(\mu|i)\notag\allowdisplaybreaks\\
&=&  \sigma_{k-1}(\mu|i)\sigma_{k-1}(\mu|i)
     -\sigma_k(\mu|i) \sigma_{k-2}(\mu|i)\notag\allowdisplaybreaks\\
&\ge&\big(1-\frac{(k-1)(n-k)}{k(n-k+1)}\big)\sigma_{k-1}^2(\mu|i)\notag\allowdisplaybreaks \\
&=&  \frac{n}{k(n-k+1)} \sigma_{k-1}^2(\mu|i)\ge0,  \label{2.1}
\end{eqnarray}
where we have used Lemma \ref{lemma1}.
Thus
\begin{eqnarray}
\frac{\partial G_k}{\partial\eta_j}\notag
&=&\frac{\partial G_k}{\partial\mu_j}\frac{\partial\mu_j}{\partial\eta_i}\notag\\
&=&(t\delta_{ij}+(1-t))\frac{\partial G_k}{\partial\mu_j}\notag\\
&=&(t\delta_{ij}+(1-t))
\frac{\sigma_{k-1}(\mu|i)\sigma_{k-1}(\mu)-\sigma_k(\mu)\sigma_{k-2}(\mu|i)}{\sigma_{k-1}^2(\mu)}\ge0,\label{2.4}
\end{eqnarray}
Now let $G_0(\eta)=- \frac{1}{\sigma_{k-1}(t\eta+(1-t)tr\eta \cdot e)}$.
It is easy to see that
\begin{eqnarray}
&&\frac{\partial G_0}{\partial\eta_j}\ge0,\label{2.99}
\end{eqnarray}
From (\ref{2.4}) and (\ref{2.99}),
we know that
$\frac{\partial G}{\partial\eta_j}
=\frac{\partial G_k}{\partial\eta_j}
+h(x)\frac{\partial G_0}{\partial\eta_j}>0$ and the operator $G$ is elliptic.

(2) Denote by $F(\eta)$ the operator
$\big({\sigma_{k-1}(t\eta+(1-t)tr\eta \cdot e)}
\big)^{\frac{1}{k-1}}$. Then $-G_0=F^{-k+1}$. Let
$\mu=t\eta+(1-t)tr\eta \cdot e$.
 By direct calculation

\begin{eqnarray}
&& \frac{\partial^2 G_0}{\partial\mu_i\partial\mu_j}  \notag\allowdisplaybreaks\\
&=&  -(-k+1)(-k)F^{-k-1}\frac{\partial F}{\partial \mu_i}\frac{\partial F}{\partial \mu_j}
     -(-k+1)\frac{\partial^2 F}{\partial\mu_i\partial\mu_j}\le 0.  \label{2.3}
\end{eqnarray}
  Note that
\begin{eqnarray}
&&\frac{\partial^2 G_0}{\partial\eta_p\partial\eta_q}  \notag\allowdisplaybreaks\\
&=& \frac{\partial^2 G_0}{\partial\mu_i\partial\mu_j}(t\delta_{pi}+(1-t)) (t\delta_{qj}+(1-t)) .\label{2.5}
\end{eqnarray}
Combining with (\ref{2.3}) and (\ref{2.5})
we obtain that the matrix
$\{\frac{\partial^2 G_0}{\partial\eta_i\partial\eta_j} \}$ is negative semi-definite
and $G_0$ is concave.
Similarly, by the concavity of $\frac{\sigma_k}{\sigma_{k-1}}$,
it is not difficulty to derive the concavity of $G_k$.
Then
$G=G_k+ h(x)G_0$ is concave.$\Box$

Now we will state  the a priori estimates for the solutions of
(\ref{eq99}), the proof of which will be given in Sect 3.1,
Sect 4.1 and Sect 4.2, respectively.
\begin{lemma}\label{estn0} Assume $f>0, -Ric_{g_0}\in\Gamma_k$. Let $u$ be a solution of (\ref{eq99}) and $U$  $\in\Gamma_{k-1}$.
Then there exists constant $C$,
 depending on
 $\alpha,
 f$, $g_0$,
such that
\begin{equation}
\sup_M| u|\le C.\label{C0bdn}
\end{equation}
\end{lemma}

\begin{lemma}\label{estn} Assume $f\ge\theta>0$. Let $u$ be a solution of (\ref{eq99}) and $U$ $\in\Gamma_{k-1}$.
Then there exists constant $C$,
 depending on
 $||\alpha||_{C^2(B_r)},
 ||f||_{C^2(B_r)}$, $\sup_{B_r} |u|$, $\theta$, $g_0$,
such that
\begin{equation}
\sup_{B_{\frac{r}{2}}}|\nabla u|^2\le C.\label{C1bdn}
\end{equation}

\end{lemma}

\begin{lemma}\label{estn2} Assume $f\ge\theta>0$. Let $u$ be a solution of (\ref{eq99}) and $U$  $\in\Gamma_{k-1}$.
Then there exists constant $C$,
 depending on
 $||\alpha||_{C^2(B_r)},
 ||f||_{C^2(B_r)}$, $\sup_{B_r} |u|$, $\sup_{B_r} |\nabla u|$, $\theta,$  $g_0$,
such that
\begin{equation}
\sup_{B_{\frac{r}{2}}} |\nabla^2 u|\le C.\label{C2bdn}
\end{equation}

\end{lemma}

\textbf{Case(B)} Now we assume $f=0$ and $\alpha<0$. Consider a
family of equations
 \begin{eqnarray}
G(U)
&=&\frac{\sigma_k\big(tU +(1-t)trU g_0\big)}
     {\sigma_{k-1}\big(tU +(1-t)trU  g_0\big)}\notag\allowdisplaybreaks\\
&=&\big((1-t)\frac{\sigma_k(e)}{\sigma_{k-1}(e)}-t\alpha\big) \exp({2u}),\label{eq99B}
\end{eqnarray}
where  $ t\in [0, 1], e=(1,\cdots, 1)$.
The following are the a priori estimates for the solutions of
(\ref{eq99B}), the proof of which will be given in Sect 3.2, Sect 4.3 and Sect 4.4.
\begin{lemma}\label{estn0B} Assume $f>0$. Let $u$ be a solution of (\ref{eq99B}) and $U\in$  $\Gamma_{k}$.
Then there exists constant $C$,
 depending on
 $\alpha$,  $g_0$,
such that
\begin{equation}
\sup_M| u|\le C.\label{C0bdnB}
\end{equation}
\end{lemma}

\begin{lemma}\label{estnB} Assume $\alpha<0$. Let $u$ be a solution of
(\ref{eq99B}) and $U$ $\in\Gamma_{k}$. Then there exists constant
$C$,
 depending on
 $||\alpha||_{C^2(B_r)}
 $, $\sup_{B_r} |u|$,  $g_0$,
such that
\begin{equation}
\sup_{B_{\frac{r}{2}}}|\nabla u|^2\le C.\label{C1bdnB}
\end{equation}

\end{lemma}

\begin{lemma}\label{estn2B} Assume $\alpha<0$. Let $u$ be a solution of (\ref{eq99B}) and $U$  $\in\Gamma_{k}$.
Then there exists constant $C$,
 depending on
 $||\alpha||_{C^2(B_r)}
 $, $\sup_{B_r} |u|$, $\sup_{B_r} |\nabla u|$, $g_0$,
such that
\begin{equation}
\sup_{B_{\frac{r}{2}}} |\nabla^2 u|\le C.\label{C2bdnB}
\end{equation}

\end{lemma}

\textbf{Case(C)}
Let $\tilde{g}=\exp(-2u)g_0$. Then
 \[A_{\tilde{g}}=\nabla^2u+du\otimes du-\frac{1}{2}|\nabla u|^2g_0+A_{g_0}.\]
Now denote
 \[W=\nabla^2u+du\otimes du-\frac{1}{2}|\nabla u|^2g_0 +A_{g_0}.\]
 We consider the following:
 \begin{eqnarray}
G(W)
&=&\frac{\sigma_k\big(W  \big)}
   {\sigma_{k-1}\big(W \big)}\notag\allowdisplaybreaks -\frac{\exp{(-2ku)}{f }
    }
    {\sigma_{k-1}\big(W \big)}\notag\allowdisplaybreaks\\
&=&-\alpha(x)\exp(-2u) .\label{eq3}
\end{eqnarray}
We will prove the following theorem, the proof of which will be
given in Sect 4.5.

\begin{lemma}\label{estp1}(Theorem\ref{pes}) Assume $f(x) >\theta>0$. Let $u$ be a solution of (\ref{eq3}) and $W$  $\in\Gamma_{k-1}$.
Then there exists constant $C$,
 depending on
 $||\alpha||_{C^2(B_r)},
 ||f||_{C^2(B_r)}$, $\sup_{B_r} |u|$, $\theta$, $g_0$,
such that
\begin{equation}
\sup_{B_{\frac{r}{2}}}|\nabla u|^2+|\nabla^2 u|\le C.\label{C1bdn}
\end{equation}

\end{lemma}

\vskip30pt

\section{$C^0$ estimates}

\subsection{Proof of Lemma\ref{estn0}}
 Consider
  \begin{eqnarray}
&&\frac{\sigma_k\big(tU +(1-t)trU g_0\big)}
  {\sigma_{k-1}\big(tU +(1-t)trU  g_0\big)}
  +t\alpha \exp({2u})\notag\allowdisplaybreaks\\
&=&\frac{\big((1-t)\sigma_k(e)+tf(x)\big)\exp({2ku})}
  {\sigma_{k-1}\big(tU +(1-t)trU  g_0\big)}\allowdisplaybreaks,\label{eq999}
\end{eqnarray}
where  $ t\in [0, 1], e=(1,\cdots, 1)$,
\[U= \nabla^2 u +\frac{1}{n-2} \triangle u g_0 +|\nabla u|^2 g_0- du\otimes du -t\frac{Ric_{g_0}}{n-2}+\frac{1-t}{n}g_0.\]
Let
 $u$ be the solution of equation (\ref{eq999}).
 Suppose the maximum point of $u$ is attained
at $x_1$.
Note that
$\frac{\sigma_k}{\sigma_{k-1}}$
is concave in $\Gamma_{k-1}$ (see\cite{HS99}).
Thus
\[\frac{\sigma_k}{\sigma_{k-1}}(A+B)
\ge \frac{\sigma_k}{\sigma_{k-1}}(A)
+\frac{\sigma_k}{\sigma_{k-1}}(B)
\ge\frac{\sigma_k}{\sigma_{k-1}}(B)\]
if $A$ is positive definite and $B\in\Gamma_{k-1}$.
Since $\nabla^2 u(x_1)$ is negative definite and $\nabla u(x_1)=0$, thus
at this point
\begin{equation*}
 \frac{\sigma_k(tB_{g_0}+(1-t)tr B_{g_0}
   g_0)}{\sigma_{k-1}(tB_{g_0}+(1-t)tr B_{g_0}
   g_0)}
\ge\frac{\sigma_k(tU+(1-t)tr U
   g_0)}{\sigma_{k-1}(tU+(1-t)tr U
   g_0)},
\end{equation*}
where
$B_{g_0}=-t\frac{Ric_{g_0}}{n-2}+\frac{1-t}{n}g_0$.
Besides,
\begin{equation*}
 \frac{1}{\sigma_{k-1}(tB_{g_0}+(1-t)tr B_{g_0}
   g_0)}
\le\frac{1}{\sigma_{k-1}(tU+(1-t)tr U
   g_0)}.
\end{equation*}
Therefore,
\begin{equation*}
 C+C \exp{(2u(x_1))} \ge C \exp{(2ku(x_1))}
\end{equation*}
and
we have
$$u(x_1)\leq C,$$
which means
$$\sup_{x \in M}u(x)\leq C.$$
Similarly, calculate at the minimal point of $u$, we obtain
$$\inf_{x \in M}u(x)\geq -C.$$
\subsection{Proof of Lemma\ref{estn0B}}
Let consider
\begin{eqnarray}
G(U)
&=&\frac{\sigma_k\big(tU +(1-t)trU g_0\big)}
     {\sigma_{k-1}\big(tU +(1-t)trU  g_0\big)}\notag\allowdisplaybreaks\\
&=&\big((1-t)\frac{\sigma_k(e)}{\sigma_{k-1}(e)}-t\alpha\big) \exp({2u}),
\end{eqnarray}
where
 $ t\in [0, 1], e=(1,\cdots, 1)$,
\[U= \nabla^2 u +\frac{1}{n-2} \triangle u g_0 +|\nabla u|^2 g_0- du\otimes du -t\frac{Ric_{g_0}}{n-2}+\frac{1-t}{n}g_0.\]
Assume $x_1$ is the maximum point of $u$.
Then

\begin{equation*}
 C\ge C \exp{(2u(x_1))}
\end{equation*}
and
we have
$$u(x_1)\leq C,$$
which means
$$\sup_{x \in M}u(x)\leq C.$$
Similarly, calculate at the minimal point of $u$, and we find
$$\inf_{x \in M}u(x)\geq -C.$$

\vskip30pt

\section{$C^1$ and $C^2$ estimates}
\subsection{Proof of Lemma \ref{estn}}
Let
\[U= \nabla^2 u +\frac{1}{n-2} \triangle u g_0 +|\nabla u|^2 g_0- du\otimes du -t\frac{Ric_{g_0}}{n-2}+\frac{1-t}{n}g_0.\]
We consider the following equations
\begin{equation}\label{102}
 \frac{\sigma_k(tU+(1-t)tr Ug_0)}{\sigma_{k-1}(tU+(1-t)tr Ug_0)}
 +t\alpha \exp{(2u)}(x)
=\frac{\big( (1-t)\sigma_k(e)+tf(x)\big)\exp{(2ku)}}{\sigma_{k-1}(tU+(1-t)tr Ug_0)},
\end{equation}
where $t \in [0, 1]$, $e=(1,\cdots, 1)$.
For the convenience of notations, we will denote
\begin{equation*}
G_k(U)=\frac{\sigma_k(tU+(1-t)trUg_0)}{\sigma_{k-1}(tU+(1-t)trUg_0)}, G_0(U)=-\frac{1}{\sigma_{k-1}(tU+(1-t)trU g_0)},
\end{equation*}
\begin{eqnarray*}
P_k(V)=\frac{\sigma_k(V)}{\sigma_{k-1}(V)}, && P_0(V)=-\frac{1}{\sigma_{k-1}(V)},\notag
\end{eqnarray*}
\begin{eqnarray*}
P(V)&=&P_k(V)+ \big( (1-t)\sigma_k(e)+tf(x)\big)\exp{(2ku)}P_0(V),
\end{eqnarray*}
where $V=tU+(1-t)trU g_0.$
We further denote by $\sigma_l^{ij}, P^{ij}$, $P_k^{ij}$, $P_0^{ij}$ $G^{ij}$, $G_k^{ij}$ and $G_0^{ij}$ the functions
$\frac{\partial \sigma_l}{\partial V_{ij}},$
$\frac{\partial P}{\partial V_{ij}}$,
$\frac{\partial P_k}{\partial V_{ij}}$,
$\frac{\partial P_0}{\partial V_{ij}}$ ,
$\frac{\partial G}{\partial U_{ij}}$,
$\frac{\partial G_k}{\partial U_{ij}}$ and
$\frac{\partial G_0}{\partial U_{ij}}$ respectively.
By direct calculation, we have
\begin{eqnarray}
&&P^{ii}V_{ii}=P_k^{ii}V_{ii}+ \big( (1-t)\sigma_k(e)+tf(x)\big)\exp{(2ku)}P_0^{ii}V_{ii}\notag\\
&=&\frac{\sigma_k^{ii}V_{ii}\sigma_{k-1}-\sigma_k\sigma_{k-1}^{ii}V_{ii}}
  {\sigma_{k-1}^2}
  + \big( (1-t)\sigma_k(e)+tf(x)\big)\exp{(2ku)}\frac{\sigma_{k-1}^{ii}V_{ii}}{\sigma_{k-1}^2}\notag\\
&=&\frac{\sigma_k\sigma_{k-1}}{\sigma_{k-1}^2}+(k-1)\big( (1-t)\sigma_k(e)+tf(x)\big)\exp{(2ku)}G_0\notag\\
&\ge&P=-t\alpha \exp{(2u)}.\label{100}
\end{eqnarray}
It is obviously that
\begin{equation}
G^{ii}U_{ii}=P^{jj}\big(t\delta_{ij}+(1-t)\big)U_{ii}=P^{jj}V_{jj}.\label{103}
\end{equation}
Moreover, we differentiate the equation (\ref{102}) and obtain
\begin{eqnarray}
&&G^{ij}U_{ijp}+ \Big(\big( (1-t)\sigma_k(e)+tf(x)\big)\exp{(2ku)}\Big)_pG_0\notag\\
&=&G_k^{ij}U_{ijp}+ \big( (1-t)\sigma_k(e)+tf(x)\big)\exp{(2ku)}G_0^{ij}U_{ijp}
  \notag\\
&&+
  \big(\big( (1-t)\sigma_k(e)+tf(x)\big)\exp{(2ku)}\big)_pG_0\notag\\
&=&-(t\alpha \exp{(2u)})_p.\label{99}
\end{eqnarray}

Set $\min_M\gamma'>0,\min_M   (\gamma''-\gamma'^2)>0$ and
 $$Q=\rho\cdot(1+ \frac{|\nabla u|^2}{2})e^{\gamma(u)}:=\rho\cdot K.$$
 Here $0\leq\rho
\leq1$ is a cutoff function depending only on $r$ such that $\rho=1$
in $B_{\frac{r}{2}}$ and $\rho=0$ outside $B_{r}$, moreover
$$|\nabla\rho|\leq\frac
{C\rho^{1/2}}{r}, \quad |\nabla^{2}\rho|\leq\frac{C}{r^{2}}.$$
 Assume that
 $\max_{M} Q=Q(\widetilde{x})$,
${U}_{ij}(\widetilde{x})$ and hence $P^{ij}(\widetilde{x})$ and $G^{ij}(\widetilde{x})$ are diagonal.
Then differentiate $K$ at the point $\widetilde{x}$.
By direct calculation, we have
\begin{eqnarray*}
K_i (\widetilde{x}) =  e^{\gamma(u)} \left((1+ \frac{u_l^2}{2})\gamma' u_i + u_l u_{li}\right),\label{guij}
\end{eqnarray*}
and
 \begin{eqnarray}
 &&K_{ij} (\widetilde{x}) \notag\\
 &=&   e^{\gamma(u)} \left((1+ \frac{u_l^2}{2})\Big((\gamma') ^2 u_i u_j
     + \gamma' u_{ij} +\gamma'' u_i u_j\Big)
     +2u_l u_{lj}\gamma' u_i
     + u_{lj} u_{li} + u_l u_{lij}\right),\label{101}
\end{eqnarray}
Note that $\widetilde{x}$ is the maximum point of $Q$, we have
\begin{equation}
0=Q_i(\widetilde{x})=\rho_i K+\rho K_i.\label{108}
\end{equation}
Moreover,
\begin{equation}
Q_{ii}=\rho_{ii} K+2\rho_iK_i+\rho K_{ii}.\label{109}
\end{equation}
If we plug (\ref{108}) into (\ref{109}), we have
\begin{equation}
Q_{ii}=\rho_{ii} K-2\rho_i\frac{\rho_i}{\rho}K+\rho K_{ii}.\label{1099}
\end{equation}
Since $G^{ij}$ is positive definite and $Q_{ij}$ is negative definite,
we find
\begin{eqnarray}
0&\ge&\rho e^{-\gamma(u)}\big( {G}^{ii} +\frac{\sum G^{pp}}{n-2}g_0^{ii}\big)Q_{ii}(\widetilde{x})\notag\\
   &\ge&\rho e^{-\gamma(u)} \big( {G}^{ii} +\frac{\sum G^{pp}}{n-2}g_0^{ii}\big) (\rho K_{ii}+\rho_{ii}K-2\rho_i\frac{\rho_i}{\rho}K) \notag\allowdisplaybreaks\\
&\ge&e^{-\gamma(u)}\big( {G}^{ii} +\frac{\sum G^{pp}}{n-2}g_0^{ii}\big)(\rho^2K_{ii}-C\rho K ).\label{1098}
\end{eqnarray}
Then inserting (\ref{101}) into (\ref{1098}) and using Ricci identity , we obtain
\begin{eqnarray}
\notag\allowdisplaybreaks\\
0&\ge&  \rho^2 \big( {G}^{ii} +\frac{\sum G^{pp}}{n-2}g_0^{ii}\big)\notag\\
&&\left(u_l u_{lii}+ (1+ \frac{u_l^2}{2})\Big(\big(\gamma'\gamma'+ \gamma''\big)u_i u_i
   + \gamma' u_{ii}  \Big)
   +2\gamma' u_l u_{li} u_i
   + u_{li} u_{li} \right)\notag\allowdisplaybreaks\\
&& -C\rho\sum_i  {G}^{ii} (|\nabla u|^2+1)
   \notag\allowdisplaybreaks\\
&\ge&  \rho^2 \big( {G}^{ii} +\frac{\sum G^{pp}}{n-2}g_0^{ii}\big)\notag\\
&&\left(u_l u_{iil}+ (1+ \frac{u_l^2}{2})\Big(\big(\gamma'\gamma'+ \gamma''\big)u_i u_i
  + \gamma' u_{ii}  \Big)
  +2\gamma' u_l u_{li} u_i
  + u_{li} u_{li} \right)\notag\allowdisplaybreaks\\
&& -C\rho\sum_i  {G}^{ii} (|\nabla u|^2+1).
   \allowdisplaybreaks\label{1097}
 \end{eqnarray}
 Moreover, recalling the definition of $U$, we obtain
 \begin{eqnarray}
0&\ge&\rho^2  u_l {G}^{ii}\Big(U_{iil}  -\big(u_p^2- u_i u_i  \big)_{l} \Big)\notag\\
&&+
  \rho^2 \gamma'  {G}^{ii} \left( U_{ii}- (  u_p^2- u_i u_i  ) \right)(1+ \frac{u_l^2}{2})  \notag\\
&&+\rho^2  {G}^{ii}    (1+ \frac{u_l^2}{2})\big(\gamma'\gamma'+ \gamma''\big)\big(u_i^2+\frac{u_p^2}{n-2}  \big)\notag\\
&&+2G^{ii}\big(\gamma'u_iu_{li}u_l+\frac{\gamma'u_pu_{lp}u_l}{n-2}\big)-C\rho\sum_i  {G}^{ii} (|\nabla u|^2+1).  \allowdisplaybreaks
\end{eqnarray}
Then, combining (\ref{100}), (\ref{103}), (\ref{99}) and (\ref{108}),
it is straightforward to show that
\begin{eqnarray}
0&\ge&\rho^2       {G}^{ii}\Big(- u_l   u_p u_{pl} + u_l  u_i u_{il} +\gamma' (-u_p^2 + u_i^2) (1+ \frac{u_l^2}{2}) \Big) \notag\\
&&
+   \rho^2 {G}^{ii}\Big(  (1+ \frac{u_l^2}{2})(\gamma'^2+\gamma'') \big(u_i^2+\frac{u_p^2}{n-2}\big)
   \Big) \notag\\
   &&-2\rho^2  {G}^{ii}   (1+ \frac{u_l^2}{2})\big(\gamma'\gamma' \big)\big(u_i^2+\frac{u_p^2}{n-2}  \big)\notag\\
 &&-C\rho^{\frac{3}{2}}\sum_iG^{ii}(|\nabla u|^3+1)
    -C\rho\sum_i  {G}^{ii} (|\nabla u|^2+1) +\rho^2 G_0(|\nabla u|^2+1)\notag\allowdisplaybreaks\notag\\
&\ge& \rho^2  {G}^{ii}\Big(\gamma' (1+ \frac{u_l^2}{2}) u^2_{p}
- \gamma' (1+ \frac{u_l^2}{2}) u^2_i  \notag\\
&& +\gamma' (-u_p^2 + u_i^2) (1+ \frac{u_l^2}{2})  + (1+ \frac{u_l^2}{2})(-\gamma'^2+\gamma'')\big( u_i^2+\frac{u_p^2}{n-2}\big)
  \Big)\notag
\\
&& -C\sum_i  {G}^{ii} (\rho^{\frac{3}{2}}|\nabla u|^3+1)  +\rho^2 G_0(|\nabla u|^2+1) \notag\allowdisplaybreaks\notag\\
&\ge& \rho^2  {G}^{ii}\Big( (1+ \frac{u_l^2}{2})(-\gamma'^2+\gamma'')\big(  \frac{u_p^2}{n-2}\big)
  \Big)\notag
\\
&& -C\sum_i  {G}^{ii} (\rho^{\frac{3}{2}}|\nabla u|^3+1)  +\rho^2 G_0(|\nabla u|^2+1) \allowdisplaybreaks,\label{105}
\end{eqnarray}
where we have used $-\gamma'^2+\gamma''>0$.
Now let us divide the proof into two cases.

(A1)$\frac{\sigma_k}{\sigma_{k-1}}\le (|\nabla u|)^{\frac{1}{k}}.$ Then
from equation (\ref{102}) we have
\[-\rho\frac{\sigma_k}{\sigma_{k-1}}-t\alpha \exp(2u)\rho =\rho{(ft+(1-t)\sigma_k(e))\exp(2ku)G_0}{}.\]
For $z\in B_r(\widetilde{x})$, we have
\begin{equation}\rho G_0(\widetilde{x})
\ge\frac{-t\sup_{B_r}\alpha \exp{(2u)}-\rho(|\nabla u|)^{\frac{1}{k}}}{\inf_{B_r}\Big[(ft+(1-t)\sigma_k(e))\exp{(2ku)}\Big]}
\ge -C(\rho(|\nabla u|)^{\frac{1}{k}}+1).\label{104}\end{equation}
Note that
$\sum_iG^{ii}\ge\frac{n-k+1}{k}$ (\cite{GZ19}),
combining (\ref{105}) and (\ref{104}), we find
\[\rho|\nabla u|^2\le C.\]

(A2)$\frac{\sigma_k}{\sigma_{k-1}}\ge (|\nabla u|)^{\frac{1}{k}}.$ By Lemme\ref{lemma1} we have
\begin{equation}|G_0|\le\big(|\nabla u|\big)^{-\frac{k-1}{k}}.\label{1096}\end{equation}
As the same as before, $\sum_iG^{ii}\ge\frac{n-k+1}{k}$ (\cite{GZ19}), (\ref{1096})and
(\ref{105}) implies
\[\rho|\nabla u|^2\le C.\]

\subsection{Proof of Lemma \ref{estn2}}

Let
\[U= \nabla^2 u +\frac{1}{n-2} \triangle u g_0 +|\nabla u|^2 g_0- du\otimes du -t\frac{Ric_{g_0}}{n-2}+\frac{1-t}{n}g_0.\]
We consider the following equation
\begin{equation}\label{999}
\frac{\sigma_k(tU+(1-t)tr Ug_0)}{\sigma_{k-1}(tU+(1-t)tr Ug_0)}
+t\alpha \exp{(2u)}(x)=\frac{\big( (1-t)\sigma_k(e)+tf(x)\big)\exp{(2ku)}}{\sigma_{k-1}(tU+(1-t)tr Ug_0)},
\end{equation}
where $t \in [0, 1]$, $e=(1,\cdots, 1)$.
Take the auxiliary function $$H(x)=\rho(\Delta u+n|\nabla
u|^2):=\rho K.$$ Here $0\leq\rho \leq1$ is a cutoff function
depending only on $r$ such that $\rho=1$ in $B_{\frac{r}{2}}$ and
$\rho=0$ outside $B_{r}$, moreover
$$|\nabla\rho|\leq\frac
{C\rho^{1/2}}{r}, \quad |\nabla^{2}\rho|\leq\frac{C}{r^{2}}.$$
Assume $x_0$ is the maximum point of $H$. Then at $x_0,$
\begin{equation}\label{97}
H_i(x_0)=\rho_i K+\rho K_i=\rho_i(u_{kk}+nu_k^2)+\rho(u_{kki}+2n u_ku_{ki})=0.
\end{equation}
Differentiating both sides of (\ref{97}), we have
\begin{equation}\label{9700}
H_{ii}(x_0)=\rho_{ii} K+\rho K_{ii}+2\rho_iK_i.
\end{equation}
If we plug (\ref{97}) back into (\ref{9700}),
we obtain
\begin{equation}\label{9701}
H_{ii}(x_0)=\rho_{ii} K+\rho K_{ii}-2\rho_i\frac{\rho_i}{\rho}K.
\end{equation}
Recall
\[U= \nabla^2 u +\frac{1}{n-2} \triangle u g_0 +|\nabla u|^2 g_0- du\otimes du -t\frac{Ric_{g_0}}{n-2}+\frac{1-t}{n}g_0.\]
We will calculate in the normal coordinates which
is centered at  $x_0$.
We further assume $U_{ij}$ and hence $\frac{\partial G}{\partial U_{ij}}$
are diagonal at the point $x_0$.
Since $U\in\Gamma_2$, we have
\[|U_{ij}|\le Ctr U, tr U>0.\]
Therefore,
\begin{equation}
|u_{ij}|\le C(\Delta u+1).\label{96}
\end{equation}
In view of (\ref{96}), we may assume \[\Delta u>C.\]
For the convenience of notations, we will denote
\begin{equation*}
G_k(U)=\frac{\sigma_k(tU+(1-t)trUg_0)}{\sigma_{k-1}(tU+(1-t)trUg_0)}, G_0(U)=-\frac{1}{\sigma_{k-1}(tU+(1-t)trU g_0)},
\end{equation*}
where $e=(1,\cdots,1)$.
We further denote by $G^{ij}$, $G_k^{ij}$ and $G_0^{ij}$ the functions
$\frac{\partial G}{\partial U_{ij}}$,
$\frac{\partial G_k}{\partial U_{ij}}$ and
$\frac{\partial G_0}{\partial U_{ij}}$ respectively.
We differentiate the equation (\ref{999}) and obtain
\begin{eqnarray}
&&G^{ij}U_{ijp}+ ( \big( (1-t)\sigma_k(e)+tf(x)\big)\exp(2ku))_pG_0\notag\\
&=&G_k^{ij}U_{ijp}+ \big( (1-t)\sigma_k(e)+tf(x)\big)\exp(2ku)G_0^{ij}U_{ijp}\notag\\
&&+ ( \big( (1-t)\sigma_k(e)+tf(x)\big)\exp(2ku))_pG_0\notag\\
&=&-(t\alpha \exp{(2u)})_p.\label{95}
\end{eqnarray}
Differentiate the equation (\ref{999}) another time and we obtain
\begin{eqnarray}
&&G^{ij,rs}U_{ijp}U_{rsp}
   +G^{ij}U_{ijpp}\notag\allowdisplaybreaks\\
&&+2(\big( (1-t)\sigma_k(e)+tf(x)\big)\exp(2ku))_pG_0^{ij}U_{ijp}\notag\\&&
   + ( \big( (1-t)\sigma_k(e)+tf(x)\big)\exp(2ku))_{pp}G_0\notag\\
&=&G_k^{ij,rs}U_{ijp}U_{rsp}
   +G_k^{ij}U_{ijpp}
   + \big( (1-t)\sigma_k(e)+tf(x)\big)\exp(2ku)G_0^{ij,rs}U_{ijp}U_{rsp}\notag\\
&&+\big( (1-t)\sigma_k(e)+tf(x)\big)\exp(2ku)G_0^{ij}U_{ijpp}\notag\allowdisplaybreaks\\
&&+2(\big( (1-t)\sigma_k(e)+tf(x)\big)\exp(2ku))_pG_0^{ij}U_{ijp}\notag\\&&
   + (\big( (1-t)\sigma_k(e)+tf(x)\big)\exp(2ku))_{pp}G_0\notag\\
&=&-(t\alpha \exp{(2u)})_{pp}.\label{94}
\end{eqnarray}
Moreover, from (3.10) in \cite{GZ19}, we have
\begin{eqnarray}
-G_0^{ij,rs}U_{ijp}U_{rsp}
\ge
-\big(1+\frac{1}{k+1}\big)G_0^{-1}G_0^{ij}G_0^{rs}U_{ijp}U_{rsp}.\label{93}
\end{eqnarray}
From the positivity of $G^{ij}$ and negativity of $H_{ij}$, we find
\begin{eqnarray}
0&\ge&\rho \big( {G}^{ii} +\frac{\sum G^{pp}}{n-2}g_0^{ii}\big)H_{ii}(x_0)\label{9300}\allowdisplaybreaks\end{eqnarray}
Then plug (\ref{9701}) into (\ref{9300}), we obtain
\begin{eqnarray}
0&\ge&\rho^2\big( {G}^{ii} +\frac{\sum G^{pp}}{n-2}g_0^{ii}\big)\big(u_{ppii}+ 2nu_p u_{pii} +2nu_{pi}u_{pi}\big)-\rho CG^{ii}\Delta u.\notag\allowdisplaybreaks\end{eqnarray}
Now using Ricci identity yields
\begin{eqnarray}
0&\ge&\rho^2\big( {G}^{ii} +\frac{\sum G^{pp}}{n-2}g_0^{ii}\big)\big(u_{iipp}+2nu_pu_{iip}+2nu_{pi}^2-C\Delta u\big)-\rho CG^{ii}\Delta u\notag\allowdisplaybreaks\end{eqnarray}
It follows from the definition of $U$ and (\ref{97}) that
\begin{eqnarray}
0&\ge&\rho^2G^{ii}\big(U_{iipp}+(u_i^2)_{pp}-\big({u_l^2}\big)_{pp}\notag\allowdisplaybreaks\\
  &&+2nu_p\big(U_{iip}+(u_i^2)_p-\big({u_l^2}\big)_p\big)+ 2nu_{pi}^2+\frac{2n}{n-2}u^2_{lp}-C\Delta u\big)
  -\rho CG^{ii}\Delta u\notag\allowdisplaybreaks\\
&\ge&\rho^2G^{ii}\big(U_{iipp}+2(-2nu_iu_{ip}u_{p}+u^2_{ip})-\big({-2nu_lu_{lp}u_p+u_{lp}^2}\big)\notag\allowdisplaybreaks\\
  &&+2nu_p\big(U_{iip}+2(u_iu_{ip})-2\big({u_lu_{lp}}\big)\big)+2nu_{pi}^2+\frac{2n}{n-2}u^2_{lp}-C\Delta u\big)
-\rho CG^{ii}\Delta u\notag\allowdisplaybreaks\\
&\ge&G^{ii}(\frac{n+2}{n-2}\rho^2(\Delta u)^2-C\rho\Delta u)
+\rho^2G^{ii}U_{iipp}+2\rho^2u_pG^{ii}U_{iip}\notag\end{eqnarray}
where we have used the assumption that $\Delta u$ is sufficiently large. Then
using the concavity of $G_k$ , (\ref{95}) and (\ref{94}) we deduce that
\begin{eqnarray}
0&\ge&G^{ii}(\frac{n+2}{n-2}\rho^2(\Delta u)^2-C\rho\Delta u)\notag\\
&&-\rho^2\Big(\big( (1-t)\sigma_k(e)+tf(x)\big)\exp(2ku)\Big)G_0^{ij,rs}U_{ijp}U_{rsp}
   \notag\allowdisplaybreaks\\
&&-2\rho^2\big( \big( (1-t)\sigma_k(e)+tf(x)\big)\exp{(2ku)}\big)_pG_0^{ij}U_{ijp}\notag\\&&
   - \rho^2\big(\big( (1-t)\sigma_k(e)+tf(x)\big)\exp{(2ku)}\big)_{pp}G_0\notag\\
&&
   -\rho^2(t\alpha \exp{(2u)})_{pp}\notag\\&&
   +2\rho^2u_p\big(-\big( \big( (1-t)\sigma_k(e)+tf(x)\big)\exp{(2ku)}\big))_pG_0
-(t\alpha \exp{(2u)})_p\big)\notag\end{eqnarray}
By use of (\ref{93}), it yields
\begin{eqnarray}
0&\ge&G^{ii}(\frac{n+2}{n-2}\rho^2(\Delta u)^2-C\rho\Delta u)\notag\\
&&-\rho^2\big( \big( (1-t)\sigma_k(e)+tf(x)\big)\exp{(2ku)}\big)\big(1+\frac{1}{k+1}\big)G_0^{-1}G_0^{ij}G_0^{rs}U_{ijp}U_{rsp}
   \notag\allowdisplaybreaks\notag\\
&&
   -2\rho^2\big( \big( (1-t)\sigma_k(e)+tf(x)\big)\exp{(2ku)}\big)_pG_0^{ij}U_{ijp}
   +C\rho^2 G_0\big(\Delta u+1\big)\notag\\
&\ge&G^{ii}(\frac{n+2}{n-2}\rho^2(\Delta u)^2-C\rho\Delta u)+C \rho G_0\big(\Delta u+1\big)\notag\\
&&-\rho^2\big( \big( (1-t)\sigma_k(e)+tf(x)\big)\exp{(2ku)}\big) \big(1+\frac{1}{k+1}\big)G_0^{-1}\notag\\
&&\Big(G_0^{ij}U_{ijp}
   \notag\allowdisplaybreaks
-\frac{\big( \big( (1-t)\sigma_k(e)+tf(x)\big)\exp{(2ku)}\big)_p}{-\big( (1-t)\sigma_k(e)+tf(x)\big) \big(1+\frac{1}{k+1}\big)G_0^{-1}}\Big)^2 \notag\\
&&-\frac{\rho^2\big(\big( (1-t)\sigma_k(e)+tf(x)\big)\exp{(2ku)}\big)_p^2}{-\big( \big( (1-t)\sigma_k(e)+tf(x)\big)\exp{(2ku)}\big) \big(1+\frac{1}{k+1}\big)G_0^{-1}}-C\notag\\
&&\notag\\
&\ge&G^{ii}(\frac{n+2}{n-2}\rho^2(\Delta u)^2-C\rho\Delta u)+C \rho^2G_0\big(\Delta u+1\big)-C.\label{920}
\end{eqnarray}
Let us divide the proof into two cases.

(A1)$\frac{\sigma_k}{\sigma_{k-1}}\le (\Delta u)^{\frac{1}{k}}$.
Then by the equation we have
\begin{equation}\rho|G_0|=\frac{\rho}{\sigma_{k-1}}\le \frac{\rho\frac{\sigma_k}{\sigma_{k-1}}+t\rho\alpha\sup_{B_r} \exp{(2u)}}{\inf_{B_r}\Big[(tf(x)+(1-t)\sigma_k(e))\exp{(2ku)}\Big]}
\le C\rho(\Delta u)^{\frac{1}{k}}.\label{921}\end{equation}
Since $\sum_i G^{ii}\ge\frac{n-k+1}{k}$,  (\ref{920}) and (\ref{921}) imply
$\rho\Delta u\le C$.

(A2) $\frac{\sigma_k}{\sigma_{k-1}}> (\Delta u)^{\frac{1}{k}}$.
Then by Lemma\ref{lemma1},
\begin{equation}|G_0|=\frac{1}{\sigma_{k-1}}\le C\big(\frac{\sigma_{k-1}}{\sigma_k}\big)^{k-1}\le(\Delta u)^{-\frac{k-1}{k}}.\label{910}\end{equation}
Now we also derive $\rho\Delta u\le C$
from (\ref{920}) and (\ref{910}).

\vskip30pt

\subsection{Proof of Lemma \ref{estnB}}
Let
\[U= \nabla^2 u +\frac{1}{n-2} \triangle u g_0 +|\nabla u|^2 g_0- du\otimes du -t\frac{Ric_{g_0}}{n-2}+\frac{1-t}{n}g_0.\]
We consider the following equation
\begin{equation}\label{132}
\frac{\sigma_k(tU+(1-t)tr Ug_0)}{\sigma_{k-1}(tU+(1-t)tr Ug_0)}
=\Big((1-t)\frac{\sigma_k(e)}{\sigma_{k-1}(e)}-t\alpha \Big)\exp{(2u)},
\end{equation}
where $t \in [0, 1]$, $e=(1,\cdots, 1)$.
For the convenience of notations, we will denote
\begin{eqnarray*}
P_k(V)=\frac{\sigma_k(V)}{\sigma_{k-1}(V)},\ G(U)=\frac{\sigma_k(tU+(1-t)trUg_0)}{\sigma_{k-1}(tU+(1-t)trUg_0)},\notag\\
\end{eqnarray*}
where $V=tU+(1-t)trU g_0.$
We further denote by $\sigma_l^{ij}$, $P_k^{ij}$, $G^{ij}$ the functions
$\frac{\partial \sigma_l}{\partial V_{ij}},$
$\frac{\partial P_k}{\partial V_{ij}}$
and
$\frac{\partial G}{\partial U_{ij}}$respectively.
By direct calculation
\begin{eqnarray}
G^{ii}U_{ii}&=&P_k^{jj}\big(t\delta_{ij}+(1-t)\big)U_{ii}\notag\\
&=& P_k^{ii}V_{ii} \notag\\
&=&\frac{\sigma_k^{ii}V_{ii}\sigma_{k-1}-\sigma_k\sigma_{k-1}^{ii}V_{ii}}
{\sigma_{k-1}^2}
\notag\\
&=&\frac{\sigma_k\sigma_{k-1}}{\sigma_{k-1}^2}\notag\\
&\ge&P_k=\Big((1-t)\frac{\sigma_k(e)}{\sigma_{k-1}(e)}-t\alpha \Big)\exp{(2u)}.\label{130}
\end{eqnarray}

Moreover, we differentiate the equation (\ref{132}) and obtain
\begin{eqnarray}
&&G^{ij}U_{ijp}
\notag\notag\\
&=&\big(\big((1-t)\frac{\sigma_k(e)}{\sigma_{k-1}(e)}-t\alpha \big)\exp{(2u)}\big)_p.\label{129}
\end{eqnarray}

Set $\min_M \gamma'>0,\min_M   (\gamma''-\gamma'^2)>0$, and
 $$Q=\rho\cdot(1+ \frac{|\nabla u|^2}{2})e^{\gamma(u)}:=\rho\cdot K.$$
 Here $0\leq\rho
\leq1$ is a cutoff function depending only on $r$ such that $\rho=1$
in $B_{\frac{r}{2}}$ and $\rho=0$ outside $B_{r}$, moreover
$$|\nabla\rho|\leq\frac {C\rho^{1/2}}{r}, \quad
|\nabla^{2}\rho|\leq\frac{C}{r^{2}}.$$
 Assume that
 $\max_{M} Q=Q(\widetilde{x})$,
${U}_{ij}(\widetilde{x})$ and hence $P^{ij}(\widetilde{x})$ and
$G^{ij}(\widetilde{x})$ are diagonal. Then differentiating $K$ at
the point $\widetilde{x}$, we have
\begin{eqnarray*}
K_i (\widetilde{x}) =  e^{\gamma(u)} \left((1+ \frac{u_l^2}{2})\gamma' u_i + u_l u_{li}\right)\label{guij3}
\end{eqnarray*}
and
 \begin{eqnarray}
 &&K_{ij} (\widetilde{x}) \notag\\
 &=&   e^{\gamma(u)} \left((1+ \frac{u_l^2}{2})\Big((\gamma') ^2 u_i u_j
   + \gamma' u_{ij} +\gamma'' u_i u_j\Big)
   +2u_l u_{lj}\gamma' u_i
   + u_{lj} u_{li} + u_l u_{lij}\right).\label{131}
\end{eqnarray}
Since $\widetilde{x}$ is the maximum point of $Q$, we have
\begin{equation}
0=Q_i(\widetilde{x})=\rho_i K+\rho K_i.\label{138}
\end{equation}
Differentiating both sides of (\ref{138}) gives
\begin{equation}
Q_{ii}=\rho_{ii} K+2\rho_iK_i+\rho K_{ii}.\label{139}
\end{equation}
Then inserting (\ref{138}) into (\ref{139})
and using the positivity of $G^{ij}$,
the negativity of $Q_{ij}$,
we deduce that
\begin{eqnarray}
0&\ge&\rho e^{-\gamma(u)}\big( {G}^{ii} +\frac{\sum G^{pp}}{n-2}g_0^{ii}\big)Q_{ii}(\widetilde{x})\notag\\
&\ge&\rho e^{-\gamma(u)} \big( {G}^{ii} +\frac{\sum G^{pp}}{n-2}g_0^{ii}\big) (\rho K_{ii}+\rho_{ii}K-2\rho_i\frac{\rho_i}{\rho}K) \notag\allowdisplaybreaks\\
&\ge&e^{-\gamma(u)}\big( {G}^{ii} +\frac{\sum G^{pp}}{n-2}g_0^{ii}\big)(\rho^2K_{ii}-C\rho K ).\label{1390}\allowdisplaybreaks\end{eqnarray}
Then we plug (\ref{131}) into (\ref{1390}) and obtain
\begin{eqnarray}
0&\ge&  \rho^2 \big( {G}^{ii} +\frac{\sum G^{pp}}{n-2}g_0^{ii}\big)\notag\\
&&\left(u_l u_{lii}+ (1+ \frac{u_l^2}{2})\Big(\big(\gamma'\gamma'+ \gamma''\big)u_i u_i
  + \gamma' u_{ii}  \Big)
  +2\gamma' u_l u_{li} u_i
  + u_{li} u_{li} \right)\notag\allowdisplaybreaks\\
&& -C\rho\sum_i  {G}^{ii} (|\nabla u|^2+1).
\end{eqnarray}
Moreover, Ricci identity gives
\begin{eqnarray}
 0\notag\allowdisplaybreaks
 &\ge&  \rho^2 \big( {G}^{ii} +\frac{\sum G^{pp}}{n-2}g_0^{ii}\big)\notag\\
&&\left(u_l u_{iil}+ (1+ \frac{u_l^2}{2})\Big(\big(\gamma'\gamma'+ \gamma''\big)u_i u_i
  + \gamma' u_{ii}  \Big)
  +2\gamma' u_l u_{li} u_i
  + u_{li} u_{li} \right)\notag\allowdisplaybreaks\\
&& -C\rho\sum_i  {G}^{ii} (|\nabla u|^2+1)
 \notag\allowdisplaybreaks.\end{eqnarray}
Using the definition of $U$,  we have
 \begin{eqnarray}
0&\ge&\rho^2  u_l {G}^{ii}\Big(U_{iil}  -\big(u_p^2- u_i u_i  \big)_{l} \Big)\notag\\
&&+
   \rho^2 \gamma'  {G}^{ii} \left( U_{ii}- (  u_p^2- u_i u_i  ) \right)(1+ \frac{u_l^2}{2})  \notag\\
&&+\rho^2  {G}^{ii}    (1+ \frac{u_l^2}{2})\big(\gamma'\gamma'+ \gamma''\big)\big(u_i^2+\frac{u_p^2}{n-2}  \big)\notag\\
&&+2G^{ii}\big(\gamma'u_iu_{li}u_l+\frac{\gamma'u_pu_{lp}u_l}{n-2}\big)\notag\\
&&-C\rho\sum_i  {G}^{ii} (|\nabla u|^2+1). \label{1039} \allowdisplaybreaks\end{eqnarray}
By substituting (\ref{130})and (\ref{129}) into (\ref{1039})
and using (\ref{138}), we obtain
\begin{eqnarray}
0&\ge&\rho^2       {G}^{ii}\Big(- u_l   u_p u_{pl} + u_l  u_i u_{il} +\gamma' (-u_p^2 + u_i^2) (1+ \frac{u_l^2}{2}) \Big) \notag\\
  &&
  +   \rho^2 {G}^{ii}\Big(  (1+ \frac{u_l^2}{2})(\gamma'^2+\gamma'') \big(u_i^2+\frac{u_p^2}{n-2}\big)
   \Big) \notag\\
     &&-2\rho^2  {G}^{ii}   (1+ \frac{u_l^2}{2})\big(\gamma'\gamma' \big)\big(u_i^2+\frac{u_p^2}{n-2}  \big)\notag\\
   && -C\sum_i  {G}^{ii} (\rho^\frac{3}{2}|\nabla u|^3+\rho|\nabla u|^2+1) \notag\allowdisplaybreaks\notag\\
&\ge& \rho^2  {G}^{ii}\Big(\gamma' (1+ \frac{u_l^2}{2}) u^2_{p}
  - \gamma' (1+ \frac{u_l^2}{2}) u^2_i  \notag\\
  && +\gamma' (-u_p^2 + u_i^2) (1+ \frac{u_l^2}{2})  + (1+ \frac{u_l^2}{2})(-\gamma'^2+\gamma'')\big( u_i^2+\frac{u_p^2}{n-2}\big)
  \Big)\notag
  \\
  && -C\sum_i  {G}^{ii} (\rho^{\frac{3}{2}}|\nabla u|^3+\rho|\nabla u|^2+1)   \notag\allowdisplaybreaks\notag\\
&\ge& \rho^2  {G}^{ii}\Big( (1+ \frac{u_l^2}{2})(-\gamma'^2+\gamma'')\big(  \frac{u_p^2}{n-2}\big)
    \Big)\notag
  \\
 & & -C\sum_i  {G}^{ii} (\rho^{\frac{3}{2}}|\nabla u|^3+1)   \allowdisplaybreaks,\label{135}
\end{eqnarray}
where we have used $-\gamma'^2+\gamma''>0$.
Thanks to
(\ref{135}),
(\ref{C1bdnB}) is proved.
\vskip30pt

\subsection{Proof of Lemma \ref{estn2B}}
Let
\[U= \nabla^2 u +\frac{1}{n-2} \triangle u g_0 +|\nabla u|^2 g_0- du\otimes du -t\frac{Ric_{g_0}}{n-2}+\frac{1-t}{n}g_0.\]
We consider the following equation
\begin{equation}\label{1299}
\frac{\sigma_k(tU+(1-t)tr Ug_0)}{\sigma_{k-1}(tU+(1-t)tr Ug_0)}
=\big((1-t)\frac{\sigma_k(e)}{\sigma_{k-1}(e)}-t\alpha\big) \exp{(2u)},
\end{equation}
where $t \in [0, 1]$, $e=(1,\cdots, 1)$.
Take the auxiliary function $$H(x)=\rho(\Delta u+n|\nabla
u|^2):=\rho K.$$
 Here $0\leq\rho
\leq1$ is a cutoff function depending only on $r$ such that $\rho=1$
in $B_{\frac{r}{2}}$ and $\rho=0$ outside $B_{r}$, moreover
$$|\nabla\rho|\leq\frac {C\rho^{1/2}}{r}, \quad
|\nabla^{2}\rho|\leq\frac{C}{r^{2}}.$$ Assume $x_0$ is the maximum
point of $H$. We will calculate in the normal coordinates which is
centered at  $x_0$. We further assume $U_{ij}$ and hence
$\frac{\partial G}{\partial U_{ij}}$ are diagonal at the point
$x_0$. Then at $x_0,$
\begin{equation}\label{397}
H_i(x_0)=\rho_i K+\rho K_i=\rho_i(u_{kk}+nu_k^2)+\rho(u_{kki}+2n u_ku_{ki})=0.
\end{equation}

Since $U\in\Gamma_2$, we have
\[|U_{ij}|\le Ctr U, tr U>0.\]
Therefore,
\begin{equation}
|u_{ij}|\le C(\Delta u+1).\label{396}
\end{equation}
In view of (\ref{396}), without loss of generality,
we may assume \[\Delta u>C.\]
For the convenience of notations, we will denote
\begin{equation*}
G(U)=\frac{\sigma_k(tU+(1-t)trUg_0)}{\sigma_{k-1}(tU+(1-t)trUg_0)},
\end{equation*}
where $e=(1,\cdots,1)$.
We further denote by $G^{ij}$ the functions
$\frac{\partial G}{\partial U_{ij}}$.
We differentiate the equation (\ref{1299}) and obtain
\begin{eqnarray}
&&G_k^{ij}U_{ijp}\notag\\
&=&\Big(\big((1-t)\frac{\sigma_k(e)}{\sigma_{k-1}(e)}-t\alpha\big) \exp{(2u)}\Big)_p.\label{395}
\end{eqnarray}
Differentiate the equation (\ref{1299}) another time and we obtain
\begin{eqnarray}
&&G^{ij,rs}U_{ijp}U_{rsp}
   +G^{ij}U_{ijpp}\notag\allowdisplaybreaks\\
&=&\Big(\big((1-t)\frac{\sigma_k(e)}{\sigma_{k-1}(e)}-t\alpha\big) \exp{(2u)}\Big)_{pp}.\label{394}
\end{eqnarray}

By direct calculation, we have
\begin{eqnarray}
0&\ge&\rho \big( {G}^{ii} +\frac{\sum G^{pp}}{n-2}g_0^{ii}\big)H_{ii}(x_0)\notag\allowdisplaybreaks\\
&=&\rho^2\big( {G}^{ii} +\frac{\sum G^{pp}}{n-2}g_0^{ii}\big)\big(u_{ppii}+ 2nu_p u_{pii} +2nu_{pi}u_{pi}\big)-\rho CG^{ii}\Delta u\notag\allowdisplaybreaks\\
&\ge&\rho^2\big( {G}^{ii} +\frac{\sum G^{pp}}{n-2}g_0^{ii}\big)\big(u_{iipp}+2nu_pu_{iip}+2nu_{pi}^2-C\Delta u\big)-\rho CG^{ii}\Delta u,\notag\allowdisplaybreaks\end{eqnarray}
where we have used Ricci identity. Moreover, from
the definition of $U$ and (\ref{397}), we obtain
\begin{eqnarray}
0&\ge&\rho^2G^{ii}\big(U_{iipp}+(u_i^2)_{pp}-\big({u_l^2}\big)_{pp}\notag\allowdisplaybreaks\\
  &&+2nu_p\big(U_{iip}+(u_i^2)_p-\big({u_l^2}\big)_p\big)+ 2nu_{pi}^2+\frac{2n}{n-2}u^2_{lp}-C\Delta u\big)
  -\rho CG^{ii}\Delta u\notag\\
&\ge&\rho^2G^{ii}\big(U_{iipp}+2(-2nu_iu_{ip}u_{p}+u^2_{ip})-\big({-2nu_lu_{lp}u_p+u_{lp}^2}\big)\notag\allowdisplaybreaks\\
  &&+2nu_p\big(U_{iip}+2(u_iu_{ip})-2\big({u_lu_{lp}}\big)\big)+2nu_{pi}^2+\frac{2n}{n-2}u^2_{lp}-C\Delta u\big)
  -\rho CG^{ii}\Delta u\allowdisplaybreaks.\label{6789}\allowdisplaybreaks\end{eqnarray}
Now substituting (\ref{395}) and (\ref{394})into (\ref{6789})
and using the concavity of $G$, we have
\begin{eqnarray}
0
&\ge&G^{ii}(\frac{n+2}{n-2}\rho^2(\Delta u)^2-C\rho\Delta u)
   +\rho^2G^{ii}U_{iipp}+2\rho^2u_pG^{ii}U_{iip}\notag\\
&\ge&G^{ii}(\frac{n+2}{n-2}\rho^2(\Delta u)^2-C\rho\Delta u).\label{92}
\end{eqnarray}
Now we also derive $\rho\Delta u\le C$
from (\ref{92}).

\subsection{Proof of Lemma \ref{estp1}}

Let
\begin{equation}
W=\nabla^2u+du\otimes du-\frac{1}{2}|\nabla u|^2g_0+A_{g_0}.
\end{equation}
We consider the following equation
\begin{equation}
\frac{\sigma_k(W )}{\sigma_{k-1}(W )}
+ \alpha\exp(-2u)=\frac{\phi(x,u)}{\sigma_{k-1}(W)},
\end{equation}
where $t \in [0, 1]$, $$\phi(x, u)=\exp{(-2ku)f}.$$
Take the auxiliary function $$H(x)=\rho\cdot(\Delta u+|\nabla
u|^2):=\rho \cdot K,$$ where $0\leq\rho \leq1$ is a cutoff function
depending only on $r$ such that $\rho=1$ in $B_{\frac{r}{2}}$ and
$\rho=0$ outside $B_{r}$, moreover $$|\nabla\rho|\leq\frac
{C\rho^{1/2}}{r}, \quad |\nabla^{2}\rho|\leq\frac{C}{r^{2}}.$$
Assume $x_0$ is the maximum point of $H$. Then at $x_0,$
\begin{equation}\label{4.0}
\rho_i K+\rho K_i=0.
\end{equation}
We will calculate in the normal coordinates which
is centered at  $x_0$.
We further assume $W_{ij}$ and hence $\frac{\partial G}{\partial W_{ij}}$
are diagonal at the point $x_0$.
Since $W\in\Gamma_2$, we have
\[|W_{ij}|\le Ctr W, tr W>0.\]
Therefore,
\begin{equation}
|u_{ij}|+|\nabla u|^2\le C(\Delta u+1).\label{4.01}
\end{equation}
In view of (\ref{4.01}), we may assume
\[\Delta u>C.\]
For the convenience of notations, we will denote
\begin{equation*}
G_k(W)=\frac{\sigma_k(W)}{\sigma_{k-1}(W)}, G_0(W)=-\frac{1}{\sigma_{k-1}(W)}.
\end{equation*}
We further denote by $G^{ij}$, $G_k^{ij}$ and $G_0^{ij}$ the functions
$\frac{\partial G}{\partial W_{ij}}$,
$\frac{\partial G_k}{\partial W_{ij}}$ and
$\frac{\partial G_0}{\partial W_{ij}}$ respectively.
We differentiate the equation (\ref{eq3}) and obtain
\begin{eqnarray}
&&G^{ij}W_{ijp}+ \phi(x,u)_pG_0\notag\\
&=&G_k^{ij}W_{ijp}+ \phi(x,u)G_0^{ij}W_{ijp}+ \phi(x,u)_pG_0\notag\\
&=&-(\alpha\exp(-2u) )_p.\label{4.1}
\end{eqnarray}
Differentiate the equation another time and we obtain
\begin{eqnarray}
&&G^{ij,rs}W_{ijp}W_{rsp}
   +G^{ij}W_{ijpp}\notag\allowdisplaybreaks\\
&&+2\phi(x,u)_pG_0^{ij}W_{ijp}
   + \phi(x,u)_{pp}G_0\notag\\
&=&G_k^{ij,rs}W_{ijp}W_{rsp}
   +G_k^{ij}W_{ijpp}\notag\\
&&+ \phi(x,u)G_0^{ij,rs}W_{ijp}W_{rsp}
   +\phi(x,u)G_0^{ij}W_{ijpp}\notag\allowdisplaybreaks\\
&&+2\phi(x,u)_pG_0^{ij}W_{ijp}
   + \phi(x,u)_{pp}G_0\notag\\
&=&-(\alpha\exp(-2u) )_{pp}.\label{4.2}
\end{eqnarray}
Moreover, from (3.10) in \cite{GZ19}, we have
\begin{eqnarray}
-G_0^{ij,rs}W_{ijp}W_{rsp}
\ge
-\big(1+\frac{1}{k+1}\big)G_0^{-1}G_0^{ij}G_0^{rs}W_{ijp}W_{rsp}.\label{4.3}
\end{eqnarray}
Besides, from
$G^{ii}\ge\frac{n-k+1}{k}$ (see \cite{GZ19})
we have
\begin{equation}
-C\ge -CG^{ii}.\label{4.6}
\end{equation}
In view of (\ref{4.01}), we may assume $\Delta u>C$.
Using (\ref{4.0}), we have
\begin{eqnarray}
0&\ge&\rho G^{ii}H_{ii}(x_0)\notag\allowdisplaybreaks\\
&=&\rho G^{ii}(\rho_{ii}K+2\rho_i K_i +\rho K_{ii})\notag\\
&\ge&\rho^2 G^{ii}\big(u_{ppii}+ 2u_p u_{pii} +2u_{pi}u_{pi}\big)-C\rho G^{ii}\Delta. u\notag\allowdisplaybreaks\end{eqnarray}
Then by use of
Ricci identity and the definition of $W$ we obtain
\begin{eqnarray}
0&\ge&\rho^2 G^{ii}\big(u_{iipp}+2u_pu_{iip}+2u_{pi}^2-C\Delta u\big)-C\rho G^{ii}\Delta u\notag\allowdisplaybreaks\\
&\ge&\rho^2 G^{ii}\big(W_{iipp}-(u_i^2)_{pp}+\big(\frac{u_l^2}{2}\big)_{pp}
+2u_p\big(W_{iip}-(u_i^2)_p
\notag\allowdisplaybreaks\\
&&+\big(\frac{u_l^2}{2}\big)_p\big)+ 2u_{pi}^2-C\Delta u\big)-C\rho G^{ii}\Delta u . \notag\allowdisplaybreaks\end{eqnarray}
From the concavity of $G_k$, (\ref{4.1}) and (\ref{4.2}) we deduce
\begin{eqnarray}
0&\ge&\rho^2 G^{ii}\big(W_{iipp}-2(-2u_iu_{ip}u_{p}+u^2_{ip})+\big({-2u_lu_{lp}u_p+u_{lp}^2}\big)\notag\allowdisplaybreaks\\
&&+2u_p\big(W_{iip}-2(u_iu_{ip})+\big({u_lu_{lp}}\big)\big)+2u_{pi}^2-C\Delta u\big) -C\rho G^{ii}\Delta u\notag\allowdisplaybreaks\\
&\ge&\rho^2 G^{ii}((\Delta u)^2)
   +\rho^2G^{ii}W_{iipp}+2\rho^2u_pG^{ii}W_{iip} -C\rho G^{ii}\Delta u\notag\\
&\ge&\rho^2 G^{ii}((\Delta u)^2)-C\rho G^{ii}\Delta u\notag\\
&&-\rho^2\phi(x,u)G_0^{ij,rs}W_{ijp}W_{rsp}
   \notag\allowdisplaybreaks\\
&&-2\rho^2\phi(x,u)_pG_0^{ij}W_{ijp}\notag\\&&
   - \rho^2\phi(x,u)_{pp}G_0
   -\rho^2(\alpha\exp(-2u) )_{pp}\notag\\
&&+2\rho^2u_p\big(-\phi(x,u)_pG_0
   -(\alpha\exp(-2u)(x))_p\big)\notag\end{eqnarray}
Moreover, using (\ref{4.3}) we have
\begin{eqnarray}
0&\ge&\rho^2 G^{ii}((\Delta u)^2)-C\rho G^{ii}\Delta u\notag\\
&&-\rho^2\phi(x,u) \big(1+\frac{1}{p+1}\big)G_0^{-1}G_0^{ij}G_0^{rs}W_{ijp}W_{rsp}
   \notag\allowdisplaybreaks\\&&
   -2\rho^2\phi(x,u)_pG_0^{ij}W_{ijp}
   \notag\\
&&+C \rho^2G_0\big(\Delta u+1\big).\notag\\
&\ge&G^{ii}(\rho^2 (\Delta u)^2-\rho\Delta u)+C\rho^2 G_0\big(\Delta u\big)\notag\\
&&-\rho^2\phi(x,u) \big(1+\frac{1}{k+1}\big)G_0^{-1}\Big(G_0^{ij}W_{ijp}
   \notag\allowdisplaybreaks
   -\frac{\phi(x,u)_p}{-\phi(x,u) \big(1+\frac{1}{k+1}\big)G_0^{-1}}\Big)^2 \notag\\
&&-\frac{\rho^2\phi(x,u)_p^2}{-\phi(x,u) \big(1+\frac{1}{k+1}\big)G_0^{-1}}\notag\\
&&\notag\\
&\ge& G^{ii}(\rho^2(\Delta u)^2-C\rho\Delta u)+C \rho^2G_0\big(\Delta u+1\big).\label{4.4}
\end{eqnarray}
Let us divide the proof into two cases.

(1)$\frac{\sigma_k}{\sigma_{k-1}}\le (\Delta u)^{\frac{1}{k}}$.
Then
\begin{equation}\rho G_0\ge\frac{\rho}{\sigma_{k-1}}\ge \frac{\rho\big(-\frac{\sigma_k}{\sigma_{k-1}}-\sup_{B_r}\alpha\exp(-2u)\big)}{\inf_{B_r}\phi(x,u(x))}
\ge- C\rho(\Delta u)^{\frac{1}{k}}.\label{4.100}\end{equation}
Note that $\sum_iG^{ii}\ge\frac{n-k+1}{k}$,
thus (\ref{4.4}) and (\ref{4.100}) imply
$\rho\Delta u\le C$.

(2) $\frac{\sigma_k}{\sigma_{k-1}}> (\Delta u)^{\frac{1}{k}}$.
Then by Lemma\ref{lemma1},
\begin{equation}|G_0|=\frac{1}{\sigma_{k-1}}\le C\big(\frac{\sigma_{k-1}}{\sigma_k}\big)^{k-1}\le(\Delta u)^{-\frac{k-1}{k}}.\label{4.5}\end{equation}
Now we also derive $\rho\Delta u\le C$
from (\ref{4.4}) and (\ref{4.5}).

\section{Proof of Theorem \ref{main2}}

\

\textbf{Case(A)} $f>0, \alpha\le 0$. Consider the equation
\begin{eqnarray}
&&         {\sigma_k(t\eta+(1-t)tr\eta \cdot e)}+t\alpha \exp{(2u)}{\sigma_{k-1}(t\eta+(1-t)tr\eta \cdot e)}\allowdisplaybreaks\notag\\
       &=&(1-t)(\exp{{2u}})^k\sigma_k(e)+t\exp{(2ku)}f(x), \label{106}
\end{eqnarray}
where  $\eta_i$ are the eigenvalues of
$g_0^{-1}\big(\nabla^2 u +\frac{1}{n-2} \triangle u g_0 +|\nabla u|^2 g_0- du\otimes du -t\frac{Ric_{g_0}}{n-2}+\frac{1-t}{n}g_0\big)$
and $\eta=(\eta_1, \cdots, \eta_n)$.

When $t=0$, (\ref{106}) becomes
\begin{eqnarray}
&&       \frac{2n-2}{n-2} \Delta u+({n}{}-1)|\nabla u|^2+1= \exp{2u}. \label{107}
\end{eqnarray}
Assume $x$ and $y$ are the maximum and minimum points of $u$ respectively. Then
by (\ref{107}),
\[1\ge \exp(2u(x))\]
and
\[1\le\exp(2u(y)).\]
Thus $u\equiv 0$.
 In other words, $u=0$ is the unique solution for (\ref{106}) at $t=0$.
 Let $F=G_k+((1-t)\sigma_k(e)+tf)\exp(2ku)G_0+t\alpha \exp{(2u)}:=G+t\alpha\exp(2u)$
 and $u_s$ be the variation of $u$ such that $u'=\phi$ at $s=0$.
 Then

 \begin{eqnarray}
 &&F'\notag\\&=&G^{ij}\phi_{ij}+\textrm{1st\ derivatives\ in \ }\phi \notag\\
 &&-\Big(-2t\alpha \exp(2u)-2k((1-t)\sigma_k(e)+tf)G_0\exp(2ku)\Big)\phi.
 \end{eqnarray}
 Thus linearized operator is invertible since $\alpha\le0$ and $G_0<0$.(See Theorem 6.14 of\cite{GT})
 Therefore the degree is nonzero.
( The Leray-Schauder degree is defined in \cite{Li89}.)
Thus after establishing the  a priori estimates Lemma\ref{estn0}, Lemma\ref{estn} and Lemma\ref{estn2}
we know that
 (\ref{106}) is uniformly elliptic.
 In particular, from
  \cite{Eva82} and \cite{Kry83},
  $u\in C^{2,\alpha}$, and the Schauder estimates
  give classical regularity.
Lastly, by homotopy-invariance
we obtain a solution at $t=1$.

\textbf{Case(B)} $f=0,\alpha<0.$

Consider the equation
\begin{eqnarray}
&&         {\frac{\sigma_k}{\sigma_{k-1}}(t\eta+(1-t)tr\eta \cdot e)}\allowdisplaybreaks\notag\\
       &=&\Big((1-t)\frac{\sigma_k(e)}{\sigma_{k-1}(e)}-t\alpha\Big)(\exp{{2u}}), \label{406}
\end{eqnarray}
where  $\eta_i$ are the eigenvalues of
$g_0^{-1}\big(\nabla^2u-du\otimes du+\frac{1}{2}|\nabla u|^2g_0 -tA_{g_0}+\frac{ 1-t}{n}g_0\big)$
and $\eta=(\eta_1, \cdots, \eta_n)$.

When $t=0$, (\ref{406}) becomes
\begin{eqnarray}
&&        \frac{2n-2}{n-2} \Delta u+({n}{}-1)|\nabla u|^2+1= \exp{2u}. \label{407}
\end{eqnarray}
Assume $x$ and $y$ are the maximum and minimum points of $u$ respectively. Then
by (\ref{407}),
\[1\ge \exp(2u(x))\]
and
\[1\le\exp(2u(y)).\]
Thus $u\equiv 0$ is the unique solution for (\ref{406}) at $t=0$.
 Besides, the linearized operator is invertible.
 Therefore the degree is nonzero.
( The Leray-Schauder degree is defined in \cite{Li89}.)
The  a priori estimates Lemma\ref{estn0B},Lemma\ref{estnB} and Lemma\ref{estn2B} imply that
 (\ref{406}) is uniformly elliptic.
 In particular, from
  \cite{Eva82} and \cite{Kry83},
  $u\in C^{2,\alpha}$, and the Schauder estimates
  give classical regularity.
Lastly, by homotopy-invariance
we obtain a solution at $t=1$.

\end{document}